\documentclass[twoside,11pt]{article}
\usepackage{jnm-preprint}
\usepackage{amsmath, amssymb}
\usepackage{graphicx,hyperref}

\newcommand{\beq}{\begin{eqnarray}}
\newcommand{\eeq}{\end{eqnarray}}
\newcommand{\be}[1]{\begin{equation}\label{#1}}
\newcommand{\ee}{\end{equation}}

\newcommand{\eps}{{\eps}}

\newcommand{\R}{{\mathbb R}}

\newcommand{\Sp}{\mathbb S}
\renewcommand{\H}{\mathrm H}
\renewcommand{\L}{\mathrm L}
\def\1{\mathbb I}
\renewcommand{\(}{\left(}
\renewcommand{\)}{\right)}
\newcommand{\nrm}[2]{\|{#1}\|_{\L^{#2}(\R^N)}}
\newcommand{\nrmC}[2]{\|{#1}\|_{\L^{#2}(\mathcal C)}}
\newcommand{\nrmR}[2]{\|{#1}\|_{\L^{#2}(\R)}}
\newcommand{\nrmRd}[2]{\|{#1}\|_{\L^{#2}(\R^d)}}
\newcommand{\inC}[1]{\int_{\mathcal C}#1\;dy}
\newcommand{\ird}[1]{\int_{\R^d}#1\;dx}
\renewcommand{\eps}{\varepsilon}

\newcommand{\nrmcnd}[2]{\|{#1}\|_{\L^{#2}(\mathcal C)}}

\renewcommand{\S}{{\mathbb S}^{d-1}}

\usepackage{color}

\newcommand{\two}[2]{\begin{array}{c}#1\cr#2\end{array}}

\author{Jean DOLBEAULT\thanks{Ceremade (UMR CNRS no.~7534), Univ.~Paris-Dauphine, Pl.~de Lattre de Tassigny, 75775 Paris Cedex~16, France. E-mail: \texttt{dolbeaul@ceremade.dauphine.fr}}\comma
Maria J. ESTEBAN\thanks{Ceremade (UMR CNRS no.~7534), Univ.~Paris-Dauphine, Pl.~de Lattre de Tassigny, 75775 Paris Cedex~16, France. E-mail: \texttt{esteban@ceremade.dauphine.fr}}
\support{This work has been supported by the ANR project NoNAP. J.D.~also acknowledges support from the ANR project CBDif-Fr.}}
\sauthor{Jean Dolbeault, and Maria J. Esteban}

\title{A scenario for symmetry breaking\\ in Caffarelli-Kohn-Nirenberg inequalities}
\stitle{Symmetry breaking in Caffarelli-Kohn-Nirenberg inequalities}

\abstract{The purpose of this paper is to explain the phenomenon of symmetry breaking for optimal functions in functional inequalities by the numerical computations of some well chosen solutions of the corresponding Euler-Lagrange equations. For many of those inequalities it was believed that the only source of symmetry breaking would be the instability of the symmetric optimizer in the class of all admissible functions. But recently, it was shown by an indirect argument that for some Caffarelli-Kohn-Nirenberg inequalities this conjecture was not true. In order to understand this new symmetry breaking mechanism we have computed the branch of minimal solutions for a simple problem. A reparametrization of this branch allows us to build a scenario for the new phenomenon of symmetry breaking. The computations have been performed using \emph{Freefem++}.}

\keywords{ground state; Schr\"odinger operator; Caffarelli-Kohn-Nirenberg inequality; radial symmetry; symmetry breaking; Roothan method; self-adaptive mesh; fixed point; bifurcation; finite element method; \emph{Freefem++}
\\[4pt] MSC 2010: 26D10; 46E35; 58E35 ; 65K10 ; 65N30
} 

\begin{document}
\maketitle
\ack{We thank Fr\'ed\'eric Hecht and Olivier Pironneau for their help for implementing our numerical scheme.}

\section{Symmetry breaking for Caffarelli-Kohn-Nirenberg inequalities}\label{Sec:Intro}

In this paper we are interested in understanding the \emph{symmetry breaking} phenomenon for the \emph{extremals} of a family of interpolation inequalities that have been established by Caffarelli, Kohn and Nirenberg in~\cite{Caffarelli-Kohn-Nirenberg-84}. More precisely, let $d\ge1$, $\theta\in(0,1]$ and define
\[
\Theta(p,d):=d\,\frac{p-2}{2\,p}\,,\;a_c:=\frac{d-2}2\,,\;b=a-a_c+\frac dp\,,\;\;\mbox{and}\;p^*(\theta,d):=\frac{2\,d}{d-2\,\theta}\,.
\]
Notice that $2^*=p^*(1,d)= \frac{2\,d}{d-2}$ if $d\ge 3$, while we set $2^*=\infty$ if $d=1$, $2$. For any $d\ge3$, we have
\[
0\le\Theta(p,d)\le\theta\le1\quad\Longleftrightarrow\quad 2\le p\le p^*(\theta,d)\le 2^*\,.
\]
We shall assume that $\theta\in(0,1]$, $a<a_c$, and $2\le p\le p^*(\theta,d)$ if either $d\ge 3$, or $d=2$ and $\theta<1$, or $d=1$ and $\theta<1/2$. Otherwise, we simply assume $2\le p<\infty$: $b\in[a,a+1]$ if $d\ge3$, $b\in(a,a+1]$ if $d=2$, and $b\in(a+\frac 12,a+1]$ if $d=1$. Under these assumptions, the \emph{Caffarelli-Kohn-Nirenberg inequality} amounts to
\be{Ineq:CKN}
\(\ird{\frac{|w|^p}{|x|^{bp}}}\)^\frac 2p\le \frac{\mathsf K_{\rm CKN}(\theta,p,\Lambda)}{|\Sp^{d-1}|^{(p-2)/p}}\(\ird{\frac{|\nabla w|^2}{|x|^{2a}}}\)^\theta\(\ird{\frac{|w|^2}{|x|^{2\,(a+1)}}}\)^{1-\theta}
\ee
with $\Lambda=(a-a_c^2)$, for all functions $w$ in the space obtained by completion of the set $\mathcal D(\R^d\setminus\{0\})$ of smooth functions with compact support contained in $\R^d\setminus\{0\}$, with respect to the norm
\[
w\mapsto\nrmRd{\,|x|^{-a}\,\nabla w\,}2^2+\nrmRd{\,|x|^{-(a+1)}\,w\,}2^2\,.
\]
We assume that $\mathsf K_{\rm CKN}(\theta,p,\Lambda)$ is the best constant in the above inequality. We will denote by $\mathsf K^*_{\rm CKN}(\theta,p,\Lambda)$ the best constant when the inequality is restricted to the set of radially symmetric functions. 

\medskip According to \cite{Catrina-Wang-01}, the Caffarelli-Kohn-Nirenberg inequality on $\R^d$ can be rewritten in cylindrical variables using the Emden-Fowler transformation
\[
s=\log|x|\;,\quad\omega=\frac{x}{|x|}\in\Sp^{d-1}\,,\quad u(s,\omega)=|x|^{a_c-a}\,w(x)\,.
\]
and is then equivalent to the following Gagliardo-Nirenberg-Sobolev inequality on the cylinder $\mathcal C:=\R\times\Sp^{d-1}$, namely
\be {CKNthetaC}
\nrmC up^2\le\,\mathsf K_{\rm CKN}(\theta,p,\Lambda)\,\(\nrmC{\nabla u}2^2+\Lambda\,\nrmC u2^2\)^\theta\,\nrmC u2^{2(1-\theta)}
\ee
for any $u\in H^1(\mathcal C)$. Here we adopt the convention that the measure on $\Sp^{d-1}$ is the uniform probability measure.

The parameters $a<a_c$ and $\Lambda>0$ are in one-to-one correspondence and we have chosen to make the constants $\mathsf K_{\rm CKN}$ and $\mathsf K^*_{\rm CKN}$ depend on $\Lambda$ rather than on $a$ because in the sequel of the paper we will work on the cylinder $\mathcal C$. Also, instead of working with the parameters $a$ and $b$, in the sequel we will work with the parameters $\Lambda$ and~$p$. 

Note that $u$ is an \emph{extremal} for \eqref{CKNthetaC} if and only if it is a minimizer for the energy functional
$$ u\mapsto Q^\theta_\Lambda[u]:=\frac{\left(\nrmcnd{\nabla u}2^2+\Lambda\, \nrmcnd u2^2\right)^\theta \left(\nrmcnd u2^2\right)^{1-\theta}}{\nrmcnd up^2} \,.$$

Radial symmetry of optimal functions in \eqref{Ineq:CKN}, or \emph{symmetry} to make it short, means that there are optimal functions in \eqref{CKNthetaC} which only depend on $s$. Equivalently, this means that $\mathsf K_{\rm CKN}(\theta,\Lambda,p)=\mathsf K_{\rm CKN}^*(\theta,\Lambda,p)$. On the opposite, we shall say that there is \emph{symmetry breaking} if and only if $\mathsf K_{\rm CKN}(\theta,\Lambda,p)>\mathsf K_{\rm CKN}^*(\theta,\Lambda,p)$. Notice that on the cylinder the symmetric case of the inequality is equivalent to the one-dimensional Gagliardo-Nirenberg-Sobolev inequality
\[
\nrmR up^2\le\,\mathsf K_{\rm CKN}^*(\theta,\Lambda,p)\;\(\nrmR{\nabla u}2^2+\Lambda\,\nrmR u2^2\)^\theta\,\nrmR u2^{2(1-\theta)}
\]
for any $u\in H^1(\R)$. The optimal constant $\mathsf K_{\rm CKN}^*(\theta,\Lambda,p)$ is explicit: see \cite{DDFT}.

\medskip Symmetry breaking of course makes sense only if $d\ge2$ and we will assume it is the case from now on. Let us summarize known results. Let
\[
\Lambda_{\rm FS}(p,\theta):=4\,\frac{d-1}{p^2-4}\,\frac{(2\theta-1)\,p+2}{p+2}\quad\mbox{and}\quad\Lambda_\star(p):=\frac{(N-1)\,(6-p)}{4\,(p-2)}\,.
\]
Symmetry breaking occurs for any $\Lambda>\Lambda_{\rm FS}$ according to \cite{Felli-Schneider-03,DDFT} (also see \cite{Catrina-Wang-01} for previous results and \cite{MR2437030} if $d=2$ and $\theta=1$).  This symmetry breaking is a straightforward consequence of the fact that for $\Lambda>\Lambda_{\rm FS}$, the symmetric \emph{extremals} are  saddle points of the energy functional $Q^1_\Lambda$ in the whole space, and thus cannot be even local minima.

If $\theta=1$, from \cite{DEL2011}, we know that symmetry holds for any $\Lambda\le\Lambda_\star(p)\,$. Moreover, according to \cite{DELT09}, there is a continuous curve $p\mapsto\Lambda_{\rm s}(p)$ with $\lim_{p\to2_+}\Lambda_{\rm s}(p)=\infty$ and $\Lambda_{\rm s}(p)>a_c^2$ for any $p\in (2,2^*)$ such that symmetry holds for any $\Lambda\le\Lambda_{\rm s}$ and there is symmetry breaking if $\Lambda>\Lambda_{\rm s}$. Additionally, we have that $\lim_{p\to2^*}\Lambda_{\rm s}(p)\le a_c^2$ if $d\ge3$ and, if $d=2$, $\lim_{p\to\infty}\Lambda_{\rm s}(p)=0$ and $\lim_{p\to\infty}p^2\Lambda_{\rm s}(p)=4$. The existence of this function $\Lambda_{\rm s}$ has been proven by an indirect way, and it is not explicitly known. It has been a long-standing question to decide whether the curves $p\to \Lambda_{\rm s}(p)$ and the curve
$p\to \Lambda_{\rm FS}(p, 1)$ coincide or not. This is still an open question. Let us notice that for all $p\in(2,2^*)$, $\Lambda_*(p)\le \Lambda_s(p)\le\Lambda_{\rm FS}(p,1)$ and that the difference $\Lambda_{\rm FS}(p,1)-\Lambda_*(p)$ is small, and actually smaller and smaller, when the dimension $d$ increases. For more details see \cite{DEL2011}.

According to \cite{1005}, existence of an optimal function is granted for any $\theta\in(\Theta(p,d),1)$, but for $\theta=\Theta(p,d)$, only if $\mathsf K_{\rm CKN}(\theta,\Lambda,p)>\mathsf K_{\rm GN}(p)$,
 where $\mathsf K_{\rm GN}(p)$ is the optimal constant in the Gagliardo-Nirenberg-Sobolev inequality
\[\label{ineqGN}
\nrm up^2\le\frac{\mathsf K_{\rm GN}(\theta,p,\Lambda)}{|\Sp^{d-1}|^{(p-2)/p}}\,\nrm{\nabla u}2^{2\,\Theta(p,d)}\,\nrm u2^{2\,(1-\Theta(p,d))}\quad\forall\; u\in\H^1(\R^d)\,.
\]
 A sufficient condition can be deduced, by comparison with radially symmetric functions, namely $\mathsf K_{\rm CKN}^*(\Theta(p,d),\Lambda,p)>\mathsf K_{\rm GN}(p)$, which can be rephrased in terms of~$\Lambda$ as $\Lambda<\Lambda_{\rm GN}^*(p)$ for some non-explicit (but easy to compute numerically) function $p\mapsto\Lambda_{\rm GN}^*(p)$. When $\theta=\Theta(p,d)$ and $\Lambda>\Lambda_{\rm GN}^*(p)$, extremal functions (if they exist) cannot be radially symmetric and in the asymptotic regime $p\to2_+$, this condition is weaker than $\Lambda>\Lambda_{\rm FS}(p,\Theta(p,d))$. One can indeed prove that $\lim_{p\to2_+}\Lambda_{\rm FS}(p,\Theta(p,d))>\lim_{p\to2_+}\Lambda_{\rm GN}^*(p)$. Actually, for every $\theta$ in an interval $(0, \theta_d)$, $\theta_d<1$, and for $p\in(2,2+\eps)$ for some $\eps>0$, sufficiently small, one can prove that $\Lambda_{\rm FS}(p,\Theta(p,d))>\Lambda_{\rm GN}^*(p)$; see \cite{springerlink:10.1007/s00526-011-0394-y} for more detailed statements. Hence, for $\theta\in(\Theta(p,d),1)$, close enough to $\Theta(p,d)$ and $p-2>0$, small, optimal functions exist and are not radially symmetric if $\Lambda>\Lambda_{\rm GN}^*(p)$, which is again a less restrictive condition than $\Lambda>\Lambda_{\rm FS}(p,\theta)$. See \cite{springerlink:10.1007/s00526-011-0394-y} for proofs and \cite{Oslo} for a more complete overview of known results.

\medskip The above results show that for some values of $\theta$ and $p$ there is \emph{symmetry breaking} outside the zone, parametrized by $\Lambda$, where the \emph{radial extremals} are unstable. So, in those zones the asymmetric \emph{extremals} are apart from the \emph{radial extremals} and symmetry breaking does not appear as a consequence of an instability mechanism. What is then the mechanism which makes the \emph{extremals} lose their symmetry ? The goal of this paper is to understand what is going on. A plausible scenario is provided by the numerical computations that we present in this paper. They have been done using \emph{Freefem++}.

The paper is organized as follows. In Section \ref{Sec:setting} we expose the theoretical setup of our numerical computations. In Section \ref{Sec:algorithm} we describe the algorithm. Section \ref{Sec:numerical-results} is devoted to the numerical results and their consequences. 

Our numerical method takes full advantage of the theoretical setup and provides a scenario which accounts for all known results, including the existence of non-symmetric extremal functions in ranges of the parameters for which symmetric critical points are locally stable. Although we cannot be sure that computed solutions are global optimal functions, we are able to present a convincing explanation of how symmetry breaking occurs.

\section{Theoretical setup and reparametrization of the problem for {\texorpdfstring{$\theta<1$}{theta<1}} with the problem corresponding to {\texorpdfstring{$\theta=1$}{theta=1}}}\label{Sec:setting}

Let us start this section with the case $\theta=1$ and consider the solutions to
\be{eqmu}
-\Delta u +\mu\,u= u^{p-1}\quad\mbox{in}\quad\mathcal C\,.
\ee
Any solution $u$ of \eqref{eqmu} is a critical point of $Q^1_\mu$ with critical value $Q^1_\mu[u]=\nrmC up^{p-2}$. Up to multiplication by a constant, it is also a solution to
\be{Eqn1}
-\Delta u+\mu\,u-\frac{\kappa}{\nrmcnd up^{p-2}}\,u^{p-1}=0\quad\mbox{in}\quad\mathcal C\,,
\ee
where
\[
\mu=\frac{\nrmcnd{\nabla u}2^2-\kappa\, \nrmcnd up^2}{\nrmcnd u2^2}\quad\Longleftrightarrow\quad\kappa=\frac{\nrmcnd{\nabla u}2^2+\mu\,\nrmcnd u2^2}{\nrmcnd up^2}=Q^1_\mu[u]\,.
\]
Since $Q^1_\mu[\lambda\,u]=Q^1_\mu[u]$ for any $\lambda\in\R\setminus\{0\}$, \eqref{eqmu} and \eqref{Eqn1} are equivalent. For simplicity, we will consider primarily the solutions of \eqref{eqmu}.

Let us denote by $u_{\mu,*}$ the unique positive symmetric solution of \eqref{eqmu} which achieves its maximum at $s=0$. We know from previous papers (see for instance \cite{Catrina-Wang-01} and references therein) that the positive solution is uniquely defined up to translations. As a consequence it is a minimizer of $Q^1_\mu$ among symmetric functions and $Q^1_\mu[u_{\mu,*}]=\nrmC{u_{\mu,*}}p^{p-2}=1/\mathsf K_{\rm CKN}^*(1,\mu,p)$. 

Our first goal is to study the bifurcation of non-symmetric solutions of \eqref{eqmu} from the branch $(u_{\mu,*})_\mu$ of symmetric ones. Let $f_1$ be an eigenfunction of the Laplace-Beltrami operator on the sphere $S^{d-1}$ corresponding to the eigenvalue $d-1$ and consider the Schr\"odinger operator $\mathcal H:=-\frac{d^2}{ds^2}+\mu+d-1-(p-1)\,u_{\mu,*}^{p-2}$ whose lowest eigenvalue is given by $\lambda_1(\mu)=d-1+\mu-\frac 14\,\mu\,p^2$ (see \cite{Landau-Lifschitz-67}, p.~74, for more details). As in \cite{Felli-Schneider-03}, let $\mu_{\rm FS}$ be such that $\lambda_1(\mu_{\rm FS})=0$, that is
\[
\mu_{\rm FS}=4\,\frac{d-1}{p^2-4}\,.
\]
We look for a local minimizer $u_\mu$ of $Q^1_\mu$ by expanding $u_\mu=u_{\mu,*}+\eps\,\varphi\,f_1+o(\eps)$ in terms of $\eps$, for $\mu$ in a neighborhood of $(\mu_{\rm FS})_+$. We find that $Q^1_\mu[u_{\mu,*}+\eps\,f_1\,\varphi]\sim\eps^2\,(\varphi,\mathcal H\,\varphi)_{L^2(\mathcal C)}$ as $\eps\to0$. The problem will studied with more details in~\cite{DE2012}.

It is widely believed that for $\theta=1$ the \emph{extremals} for the Caffarelli-Kohn-Nirenberg inequalities \eqref{CKNthetaC} are either the symmetric solutions $u_{\mu,*}$ for $\mu\le \mu_{\rm FS}$ or the solutions belonging to the branch which bifurcates from $u_{\rm FS}:=u_{\mu_{\rm FS},*}$. Numerically, we are going to see that this is a convincing scenario.

\medskip For $\theta<1$, the Euler-Lagrange equation for the critical points of $Q^\theta_\Lambda$ can be written as
\be{Eqntheta}
-\Delta v+\frac1\theta\Big[(1-\theta)\,t[v] +\Lambda\Big]\,v-\frac{\kappa}{\theta\,\nrmC vp^{p-2}}\,v^{p-1}=0\quad\mbox{in}\quad\mathcal C\,,
\ee
\[
\mbox{where}\quad t[v]:=\frac{\nrmcnd{\nabla v}2^2}{\nrmcnd v2^2}\quad\mbox{and}\quad \kappa=\frac{\nrmcnd{\nabla v}2^2+\Lambda\,\nrmcnd v2^2}{\nrmcnd vp^2}=Q_\Lambda^1[v]\,.
\]
Any solution $u$ of \eqref{Eqn1} is also a solution of \eqref{Eqntheta} if $\Lambda=\theta\,\mu-(1-\theta)\,t[u]$. Symmetric solutions $u_{\mu,*}$ give rise to a symmetric branch of solutions $\mu\mapsto v_{\Lambda_*^\theta(\mu),*}=u_{\mu,*}$ of solutions for \eqref{Eqntheta}, where
\[
\Lambda_*^\theta(\mu):=\theta\,\mu-(1-\theta)\,t[u_{\mu,*}]\,.
\]
A branch $\mu\mapsto u_\mu$ of solutions for \eqref{eqmu} indexed by $\mu$, normalized by the condition $\nrmC{u_\mu}p^{p-2}=\kappa$ can be seen as a branch of solutions of \eqref{Eqn1} and also provides a branch $\mu\mapsto v_{\Lambda^\theta(\mu)}=u_\mu$ of solutions for \eqref{Eqntheta} with
\[
\Lambda^\theta(\mu):=\theta\,\mu-(1-\theta)\,t[u_\mu]\,.
\]
If we can prove that $\mu\mapsto u_\mu$ bifurcates from $\mu\mapsto u_{\mu,*}$ at $\mu=\mu_{\rm FS}$, then we also have found a branch $\Lambda\mapsto v_\Lambda$ of solutions of  \eqref{Eqntheta} which bifurcates from $\Lambda\mapsto v_{\Lambda,*}$~at
\[
\Lambda_{\rm FS}^\theta:=\Lambda^\theta(\mu_{\rm FS})=\theta\,\mu_{\rm FS}-(1-\theta)\,t[u_{\rm FS}]\,
\]
as has already been noticed~in~\cite{DDFT}.

So, from the branch of solutions to \eqref{Eqn1} which bifurcates from $u_{\rm FS}$ we will construct a branch of solutions to \eqref{Eqntheta} that contains candidates to be \emph{extremals} for the Caffarelli-Kohn-Nirenberg inequalities for $\theta<1$. Of course, nothing guarantees that the \emph{extremals} for the inequalities \eqref{CKNthetaC} with $\theta<1$ lie in this branch but, as we shall see, such a scenario accounts for all theoretical results which have been obtained up to now. In this paper we numerically compute the first branch bifurcating from the symmetric branch at $\mu_{\rm FS}$, then transform it to a branch for \eqref{Eqntheta} and study the value of $J^\theta(\mu):=Q_{\Lambda(\mu)}^\theta[u_\mu]$ in terms of $\Lambda^\theta(\mu)$.

\medskip An important ingredient in our algorithm is related to the fact that minimizing the first eigenvalues of Schr\"odinger operators $-\Delta -V$ under some integral constraint on $V$ is equivalent to solving \eqref{eqmu}. One can indeed prove that
\[
\mu=-\inf_{\nrmcnd Vq=1}\quad \inf_{u\in H^1(\mathcal C)\setminus\{0\}}\frac{\inC{|\nabla u|^2}-\inC{V\,|u|^2}}{\inC{|u|^2}}
\]
with $q=\frac p{p-2}$ has a minimizer, and that the operator $-\Delta-V$ has a lowest negative eigenvalue, $-\mu$, with associated eigenfunction $u_\mu$. Up to multiplication by a constant, $V=u_\mu^{p-2}/\nrmC up^{p-2}$, so that $u_\mu$ solves \eqref{eqmu}; see \cite{DEL2011} for details, and also~\cite{MR0121101,Lieb-Thirring76}.

\section{The algorithm}\label{Sec:algorithm}

Based on the theoretical setup of the previous section, we are now ready to introduce the algorithm used for the computation of the branch of non-symmetric solutions of~\eqref{CKNthetaC} which bifurcates from $\mu_{\rm FS}$ for $\theta=1$.

\medskip\noindent\emph{1) Initialization: obtaining one point on the non-symmetric branch.} The function $u_{\mu,*}$ is a saddle point of $Q^1_\mu$ for any \hbox{$\mu>\mu_{\rm FS}$}. Choose then $\mu_0$ larger than $\mu_{\rm FS}$, but not too large, $\eps>0$ small and $w$ a direction of descent for $Q^1_{\mu_0}$ at $u_{\mu_0,*}$. Starting from $u_{\mu_0,*}+\eps\,w$, we use the conjugate gradient algorithm to decrease the energy and search for a quasi-local minimum of $Q^1_{\mu_0}$. The limit $u_{\mu_0}$ solves \eqref{Eqn1} for some $\kappa_0=Q^1_{\mu_0}[u_{\mu_0}]$. Since its energy is lower than $Q^1_{\mu_0}[u_{{\mu_0},*}]$, it is certainly a non-symmetric critical point.

\medskip\noindent\emph{2) A fixed point method.} Critical points of $Q^1_\mu$ can be characterized as fixed points of an \emph{ad hoc} map as follows.
\begin{enumerate}
\item Choose $\kappa>0$, $p\in(2,2^*)$, $q=\frac p{p-2}$ and start with some potential $V_0$ normalized by the condition: $\nrmC{V_0}q=1$.
\item For any $i\ge1$, define
\[
\lambda_i(\kappa):=\inf_{\two{u\in\H^1(\mathcal C)}{\nrmC u2=1}}\(\inC{|\nabla u|^2}-\kappa\,\inC{V_{i-1}\,|u|^2}\)\,,
\]
and get a minimizer $u_i\in\H^1(\mathcal C)$ such that $\nrmC{u_i}2=1$.
\item Define $V_i:=|u_i|^{p-2}/\nrmC{u_i}p^{p-2}$ and iterate, by computing $\lambda_{i+1}(\kappa)$ from Step~2.
\end{enumerate}
The sequence $(\lambda_i)_{i\ge1}$ is monotone non-increasing. Indeed, for any $i\ge1$, we have
\begin{eqnarray*}
\lambda_{i+1}(\kappa)&=&\inf_{\two{u\in\H^1(\mathcal C)}{\nrmC u2=1}}\(\inC{|\nabla u|^2}-\kappa\,\inC{V_i\,|u|^2}\)\\
&\le&\inC{|\nabla u_i|^2}-\kappa\,\inC{V_i\,|u_i|^2}\\
&&\qquad=\inf_{\two{V\in\L^q(\mathcal C)}{\nrmC Vq=1}}\(\inC{|\nabla u_i|^2}-\kappa\,\inC{V\,|u_i|^2}\)\\
&&\qquad\le\inC{|\nabla u_i|^2}-\kappa\,\inC{V_{i-1}\,|u_i|^2}=\lambda_i(\kappa)\,.
\end{eqnarray*}
The sequence $(\lambda_i)_{i\ge1}$ is bounded from below, as an easy consequence of H\"older's inequality:
\[
\inC{|\nabla u|^2}-\kappa\,\inC{V_i\,|u|^2}\ge\inC{|\nabla u|^2}-\kappa\,\nrmC Vq\,\nrmC up^2\,,
\]
and the r.h.s. itself is bounded by the Caffareli-Kohn-Nirenberg inequality:
\[
\inf_{\two{u\in\H^1(\mathcal C)}{\nrmC u2=1}}\(\inC{|\nabla u|^2}-\kappa\,\nrmC up^2\)=:\mu(\kappa)
\]
if we assume that $\nrmC Vq=1$. Indeed this amounts to
\[
\nrmC{\nabla u}2^2+\mu(\kappa)\,\nrmC u2^2\ge\kappa\,\nrmC up^2\quad\forall\;u\in\H^1(\mathcal C)\,,
\]
which is exactly equivalent to \eqref{Ineq:CKN} up to a reparametrization of $\mu$ in terms of $\kappa$. This scheme is converging towards a solution of 
\[
-\Delta u+\mu(\kappa)\,u=\kappa\;V\,u\;,\quad V=\kappa\,\nrmC up^{2-p}\,u^{p-2}\quad\mbox{in}\;\mathcal C\,.
\]
See \cite{DEL2011} for details and \cite{MR0121101,Lieb-Thirring76} for earlier references. Notice that if we start with the potential $V_0=u_{\mu_0}^{p-2}/\nrmC{u_{\mu_0}}p^{p-2}$ found at the end of the initialization of our algorithm, then we find that $\mu_0=\mu(\kappa_0)$ and $u_{\mu_0}$ (as well as $V_0$) is a fixed point of our algorithm, such that $\kappa_0=Q^1_{\mu_0}[u_{\mu_0}]<Q^1_{\mu_0}[u_{{\mu_0},*}]$.

The above iterative algorithm is a local version of a Roothan algorithm to compute a fixed point. We have run it using \emph{Freefem++} in a self-adaptive way, where at every step the computing mesh is based on the level lines of the previously computed function. This is important for large values of $\mu$ because solutions asymptotically tend to concentrate at some point.

\medskip\noindent\emph{3) Building the branch in terms of $\kappa$, starting from $\kappa_0$.} We adopt a perturbative approach by modifying the value of the parameter $\kappa$ and reapplying the above fixed-point algorithm.

In practice we take $\kappa=\kappa_0-\eta$ for $\eta>0$ small, $V_0=|u_{\mu_0}|^{p-2}/\nrmC{u_{\mu_0}}p^{p-2}$. We get a new critical point $u_\mu$ with $\mu=\mu(\kappa)$. If $\eta$ has been chosen small enough, we still have that $Q_\mu^1[u_\mu]<Q_\mu^1[u_{\mu,*}]$. By iterating this method as long as $\mu >\mu_{\rm FS}$, we obtain a discretized branch of numerical solutions $\mu\mapsto u_\mu$ of \eqref{Eqn1} containing~$u_{\mu_0}$. Numerically, we check that $Q_\mu^1[u_\mu]<Q_\mu^1[u_{\mu,*}]$, hence proving that $u_\mu$ is non-symmetric as long as $\mu>\mu_{\rm FS}$, and  such that $u_\mu$ converges to $u_{\rm FS}$ as $\mu$ tends to~$\mu_{\rm FS}$. For simplicity, we adopt the following convention: we extend the branch to any value of $\mu>0$ but observe that $u_\mu$ coincides with $u_{\mu,*}$ for any $\mu<\mu_{\rm FS}$. 

We can do the same in the other direction and start with $\kappa=\kappa_0+\eta$ for $\eta>0$ small, take again $V_0$ as initial potential, and then iterate, with no limitation on~$\kappa$. Again we check that $Q_\mu^1[u_\mu]<Q_\mu^1[u_{\mu,*}]$ for the discrete version of the branch corresponding to $\mu=\mu(\kappa)$, $\kappa>\kappa_0$.
 
\medskip Altogether, we have approximated a branch which bifurcates from the symmetric one at $\mu=\mu_{\rm FS}$. This branch is a very good candidate to be the branch of the global \emph{extremals} for the Caffarelli-Kohn-Nirenberg inequalities for $\theta=1$. The main reason for this belief is that if we start our algorithm for a $\mu$ close enough to $\mu_{\rm FS}$, we always hit the first (in terms of $\mu$) branch bifurcating from the symmetric one, that is, the one bifurcating from $\mu_{\rm FS}$ and observe that its energy is below the energy of corresponding symmetric solutions. On the other hand, the asymptotic value as $\mu\to\infty$ is the one predicted by Catrina and Wang in \cite{Catrina-Wang-01} for optimal functions. The estimates found in \cite{DEL2011} indicate that the set of parameters in which a different branch of optimal functions would co-exist with the branch of critical points we have computed is remarkably narrow, and close in energy with the one we have found at least for $\mu=\mu_{\rm  FS}$. The existence of another, distinct branch of non-symmetric solutions which does not bifurcate from the symmetric ones, but still has the same asymptotics as $\mu\to\infty$, is therefore very unlikely. Hence we expect that our method provides a complete answer for optimal functions and for the value of the best constant in~\eqref{Ineq:CKN} for $\theta=1$.
 
\medskip Once we have constructed a discretized branch of solutions for $\theta=1$, we use the transformation described in Section \ref{Sec:setting} to get a discretized branch of solutions for \eqref{Eqntheta}. Same comments apply as for the case $\theta=1$ (see \cite{DE2012} for the asymptotics as $\mu\to\infty$) and we expect that the computed solution is the optimal one for \eqref{Ineq:CKN} with~$\theta<1$.

\medskip Let us finally notice that by \cite{MR2001882} (also see \cite{Lin-Wang-04}), the \emph{extremals} of \eqref{CKNthetaC} enjoy some minimal symmetry properties. They depend only on $s$ and on one angle, which can be chosen as the azimuthal angle on $\S$. For any $d\ge2$, the problem we have to solve is actually two-dimensional, which greatly simplifies the computations.

\section{Numerical results}\label{Sec:numerical-results}

Solutions have been computed for various values of $p$ and $\theta\in[\Theta(p,d),1]$ in the typical case $d=5$. We adopt the convention that solutions with lowest energy are represented by plain curves while symmetric ones, when they differ, are represented by dashed curves. Darker parts of the plots correspond to minimizers for fixed $\Lambda$, at least among the branches we have computed. Hence, define the functions
\[
J^\theta(\mu):=Q_{\Lambda^\theta(\mu)}^\theta[u_\mu]\quad\mbox{and}\quad J_*^\theta(\mu):=Q_{\Lambda_*^\theta(\mu)}^\theta[u_{\mu,*}]
\]
where $u_\mu$ is the solution of \eqref{eqmu} which was obtained by the method of Section~\ref{Sec:algorithm} and $u_{\mu,*}$ is the symmetric solution of \eqref{eqmu}. For $\mu<\mu_{\rm FS}$, these two functions coincide and their value is a good candidate for determining the best constant in \eqref{Ineq:CKN}. At $\mu=\mu_{\rm FS}$, non-symmetric solutions bifurcate from symmetric ones and for $\mu>\mu_{\rm FS}$ the corresponding guess for optimal constants is given respectively by $J^\theta$ and $J_*^\theta$.

With the notations of Section~\ref{Sec:setting}, we know that  $\Lambda^{\theta=1}(\mu)=\mu$. For $\theta<1$, the parameter we are interested in is $\Lambda$ and we look for the solution of \eqref{Eqntheta} which minimizes $Q_\Lambda^\theta$. The reparametrization of Section~\ref{Sec:setting} and the bifurcation at $\mu=\mu_{\rm FS}$ suggest to still parametrize the set of solutions by $\mu$ and consider
\[
\mu\mapsto\(\Lambda^\theta(\mu),J^\theta(\mu)\)\quad\mbox{and}\quad\mu\mapsto\(\Lambda_*^\theta(\mu),J_*^\theta(\mu)\)\,.
\]
However, we will see that there is no reason why $\mu\mapsto\Lambda^\theta(\mu)$ should be monotonically increasing for $\mu>\mu_{\rm FS}$; this is indeed not the case for certain values of $d$, $p$ and $\theta$.

To illustrate these preliminary remarks, assume first that $d=5$ and $p=2.8$. For $\theta=1$, the bifurcation of non-radial solutions from the symmetric ones occurs for $\mu=\mu_{\rm FS}\approx 4.1667$. For $\theta$ ranging between $\Theta(2.8,5)\approx0.714286$ and $1$, branches can  be computed using the reparametrization of Section~\ref{Sec:setting}. See~Fig.~1. Although it is hard to see it on Fig.~1, right, the curves $\mu\mapsto\(\Lambda^\theta(\mu),J^\theta(\mu)\)$ may have a self intersection. We are now going to investigate this issue in greater details.

\medskip In the critical case $\theta=\Theta(p,d)$, the limiting value $J_\infty:=\lim_{\mu\to\infty}J^{\Theta(p,d)}(\mu)$ along the branch of non-symmetric solutions corresponds to the best constant in Gagliardo-Nirenberg inequalities:
\[
J_\infty=\mathsf k\,\inf_{u\in\H^1(\R^d)\setminus\{0\}}\frac{\ird{|\nabla u|^2}+\ird{|u|^2}}{\(\ird{|u|^p}\)^\frac 2p}
\]
where $\mathsf k:=\(s\theta^\theta\,(1-\theta)^{1-\theta}\){}_{|\theta=\Theta(p,d)}$, according to~\cite{DE2012}. When $\theta=\Theta(p,d)$, at least in cases we have computed, there are optimal functions (see Fig.~2.) if and only if
\[
\Lambda\le\Lambda_{\rm GN}:=\sup\{\Lambda_*^{\Theta(p,d)}(\mu):J_*^{\Theta(p,d)}(\mu)<J_\infty\}\,.
\]
\begin{figure}[t]\begin{center}
\includegraphics[width=6cm]{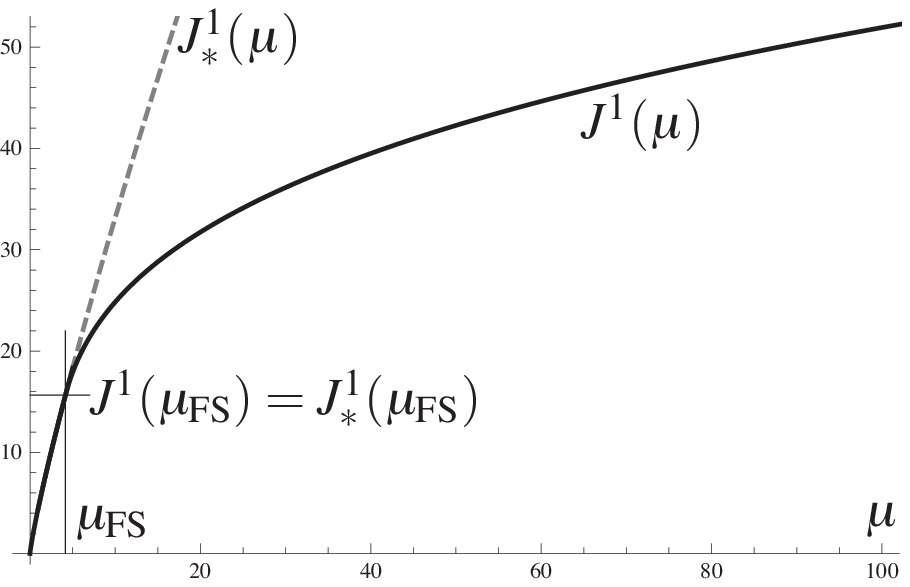}
\hspace*{1cm}
\includegraphics[width=5cm]{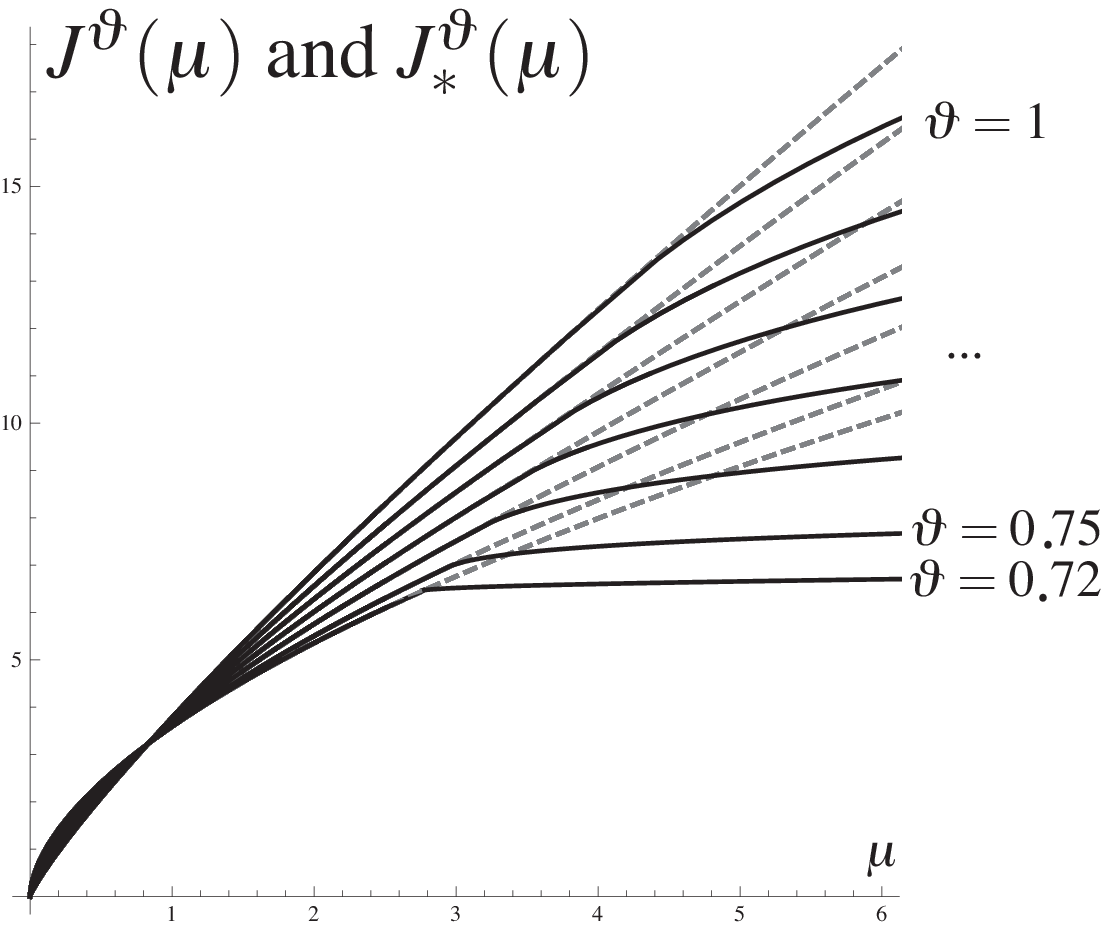}
\caption{\small\emph{{\rm Left.--} Plot of $\mu\mapsto J^1(\mu)$ (plain curve) and $\mu\mapsto J_*^1(\mu)$ (dashed curve) for \hbox{$d=5$}, $p=2.8$, $\theta=1$. The branch of non-symmetric functions bifurcates from the branch of symmetric ones for $ J^1(\mu_{\rm FS})\approx4.17$ and $Q_{\mu_{\rm FS}}^1[u_{\rm FS}]\approx15.65$. {\rm Right.--} Plots of $\mu\mapsto J^\theta(\mu)$ and $\mu\mapsto J_*^\theta(\mu)$ (dashed curve) for $d=5$, $p=2.8$, for $\theta=0.72$, $0.75$, $0.8$, $0.85$, $0.9$, $0.95$, $1$.}}
\end{center}
\vspace*{-12pt}
\end{figure}

\vspace*{-12pt}
\begin{figure}[hb]
\vspace*{-12pt}
\begin{center}
\includegraphics[width=4.5cm]{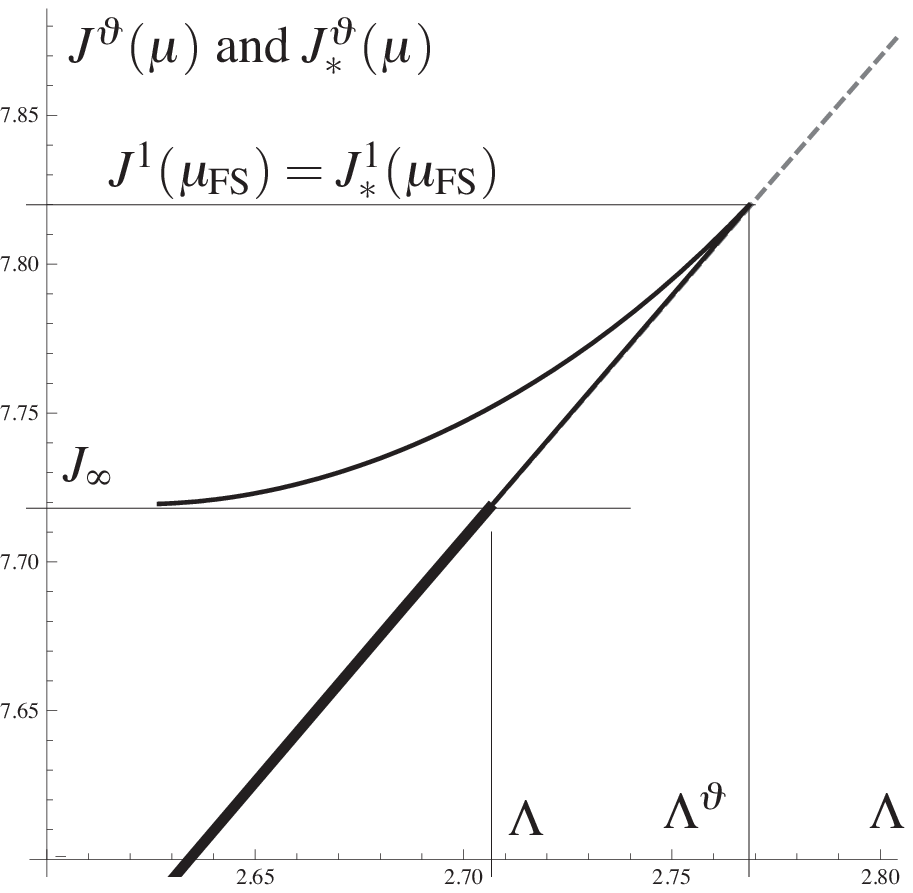}
\hspace{2cm}
\includegraphics[width=4.5cm]{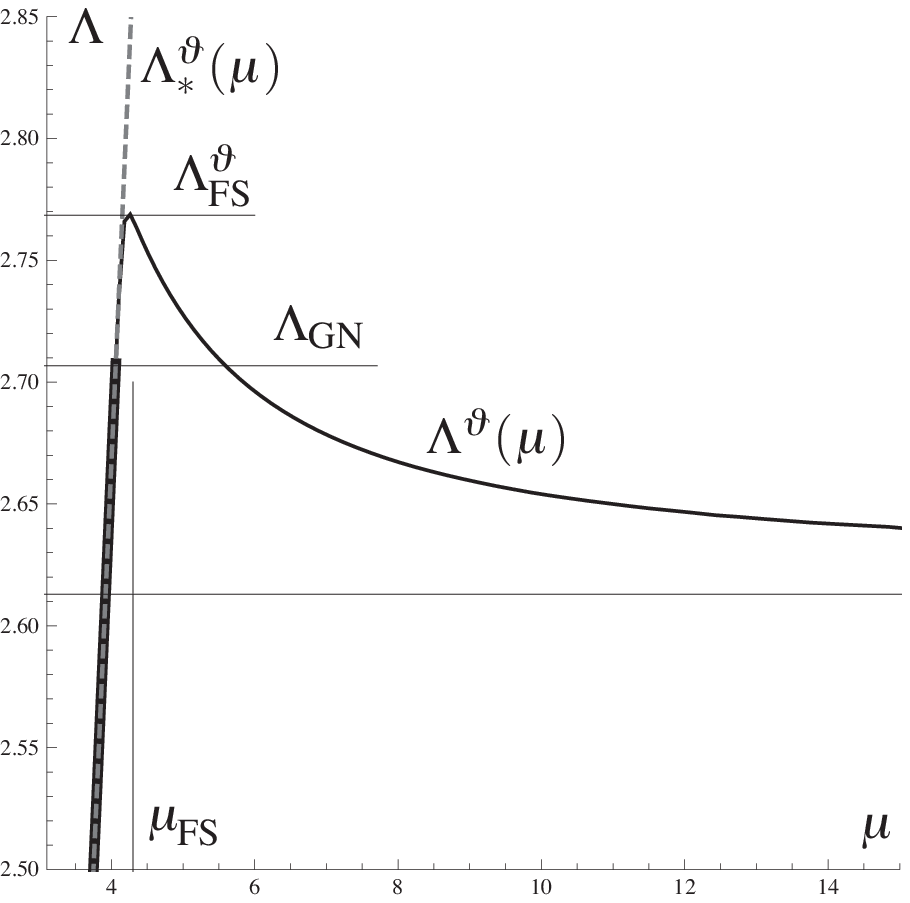}
\caption{\small\emph{Case $d=5$, $p=2.8$, $\theta=\Theta(5, 2.8)$. {\rm Left.--} Plot of $\mu\mapsto\(\Lambda^\theta(\mu),J^\theta(\mu)\)$ and (dashed curve) $\mu\mapsto\(\Lambda_*^\theta(\mu),J_*^\theta(\mu)\)$. {\rm Right.--} Reparametrization $\mu\mapsto\Lambda_*^\theta(\mu)$ (dashed curve) and $\mu\mapsto\Lambda^\theta(\mu)$: they differ for $\mu>\mu_{\rm FS}\approx4.17$.}}
\end{center}
\end{figure}

\clearpage Fix some $\theta_0\in(0,1)$. Now we consider the subcritical regime, that is when $p$ varies in the range $(2,p^*(\theta_0,d))$. When $p$ is close enough to $p^*(\theta_0,d)$, the branch stills bifurcates towards the left from the symmetric branch at $\Lambda^\theta_{\rm FS}=\Lambda^\theta(\mu_{\rm FS})$ and then, for larger values of $\mu$, turns towards the right, crosses the symmetric branch (in the plane $(\Lambda, Q^\theta_\Lambda[u])$, not in the functional space), and then stays under this symmetric branch for any larger value of $\mu$. In other words, the map $\mu\mapsto\Lambda^{\theta_0}(\mu)$ is monotone decreasing in $(\mu_{\rm FS}, \mu_0)$ for some $\mu_0>\mu_{\rm FS}$ and then monotone increasing in $(\mu_0,\infty)$. See Fig.~3.
\begin{figure}[hb]
\begin{center}
\includegraphics[width=5cm]{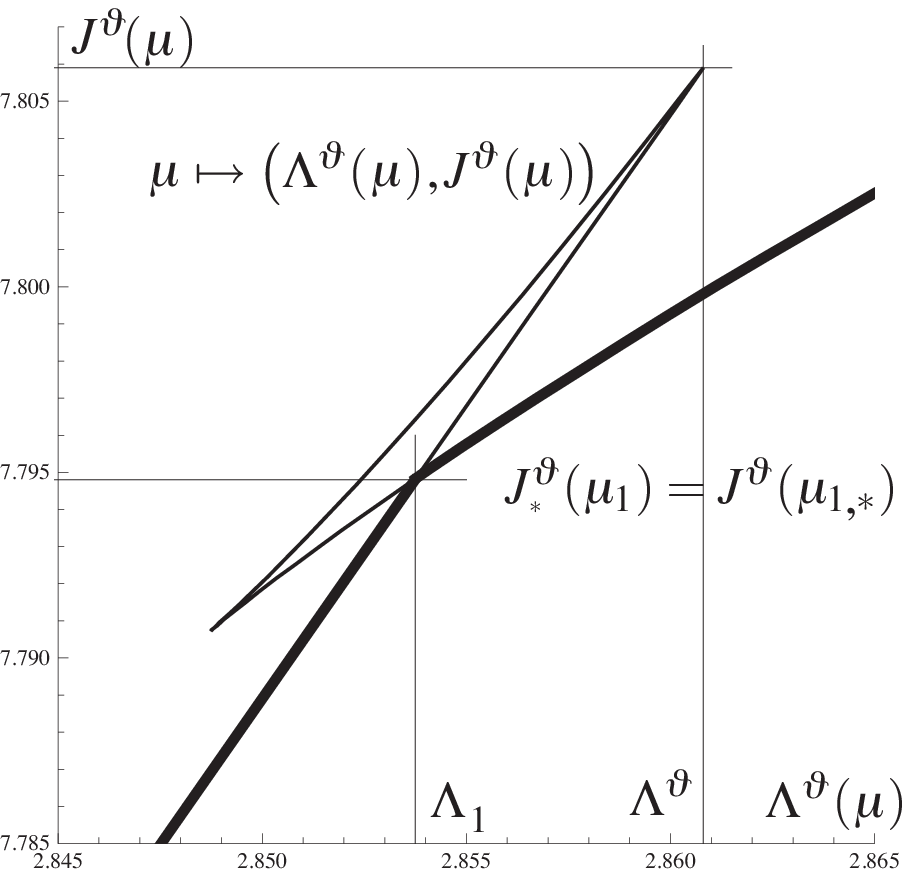}
\hspace{0.5cm}
\includegraphics[width=7cm]{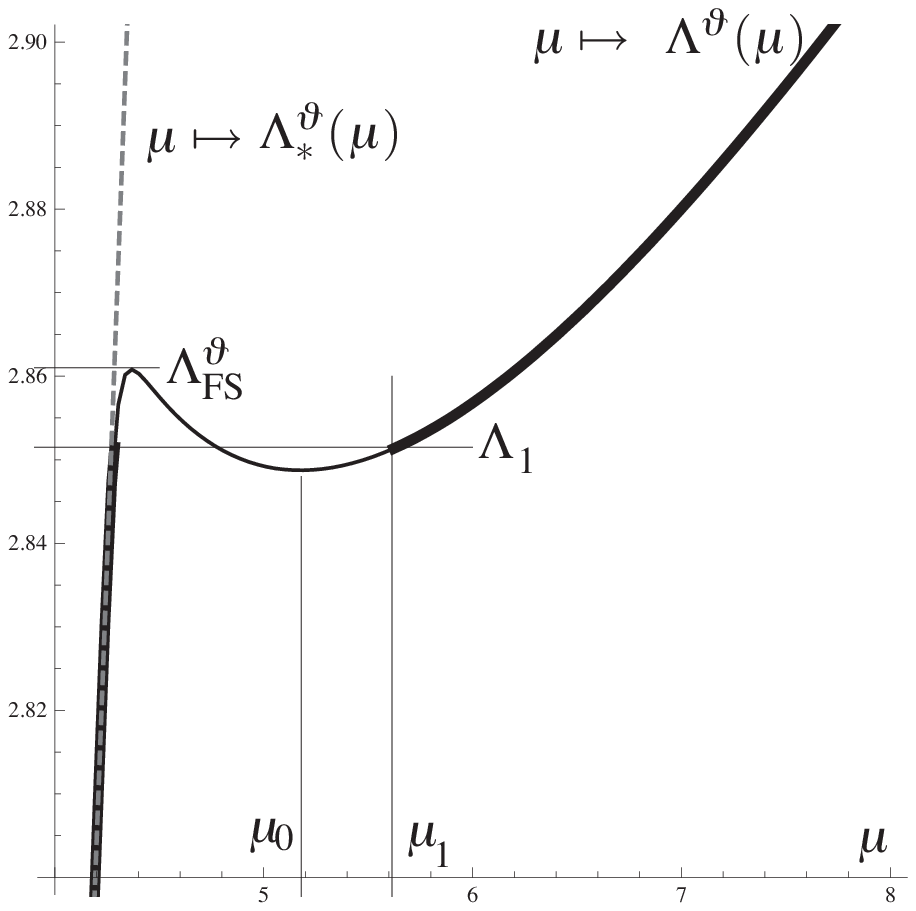}
\caption{\small\emph{Plots for $d=5$, $p=2.78$, $\theta=\Theta(5, 2.8)$. {\rm Left.--} Plot of $\mu\mapsto\(\Lambda^\theta(\mu),J^\theta(\mu)\)$. {\rm Right.--} Reparametrization $\mu\mapsto\Lambda_*^\theta(\mu)$ (dashed curve) and $\mu\mapsto\Lambda^\theta(\mu)$: the function $\Lambda_*^\theta$ is increasing for $\mu<\mu_{\rm FS}$, while the function $\Lambda^\theta$ is decreasing on $(\mu_{\rm FS}, \mu_1)$ and  increasing for $\mu>\mu_1$.}}
\end{center}
\end{figure}
This case is interesting, because when the non-symmetric branch crosses the symmetric one, if the critical points are actually the optimal functions for \eqref{Ineq:CKN}, then a symmetric \emph{extremal} and a non-symmetric one coexist. Denote by $\Lambda_1$ the corresponding value of $\Lambda$. In this case, for $\Lambda<\Lambda_1$ the \emph{extremals} of \eqref{CKNthetaC} are all symmetric. At $\Lambda=\Lambda_1$ there is coexistence of a symmetric and a non-symmetric \emph{extremal}, and for any $\Lambda > \Lambda_1$ the \emph{extremals} are non-symmetric. Let us denote by $u_\Lambda^\theta$ the minimizer of $Q^\theta_\Lambda$: the map $\Lambda\mapsto u_\Lambda^\theta$ is not continuous, since at $\Lambda =\Lambda_1$ there is a jump. 

In terms of $\mu$, there is some $\mu_1^*\in(0,\mu_{\rm FS})$ such that $\Lambda_*^\theta(\mu_1^*)=\Lambda_1$, and extremal functions for $\Lambda<\Lambda_1$ are parametrized by $\mu\mapsto(\Lambda_*^\theta(\mu),u_{\mu,*})$. There is also some $\mu_1\in(\mu_0,\infty)$ such that $\Lambda^\theta(\mu_1)=\Lambda_1$, and extremal functions for $\Lambda>\Lambda_1$ are parametrized by $\mu\mapsto(\Lambda^\theta(\mu),u_\mu)$. At $\Lambda=\Lambda_1$, $u_{\mu_1^*,*}$ and $u_{\mu_1}$ are two different optimal functions.

In Figs.~4 and 5, we show the very different shapes of the symmetric and the non-symmetric \emph{extremals} at $\Lambda =\Lambda_1$. \clearpage
\begin{figure}[ht]
\vspace*{-12pt}
\begin{center}
\includegraphics[width=5cm]{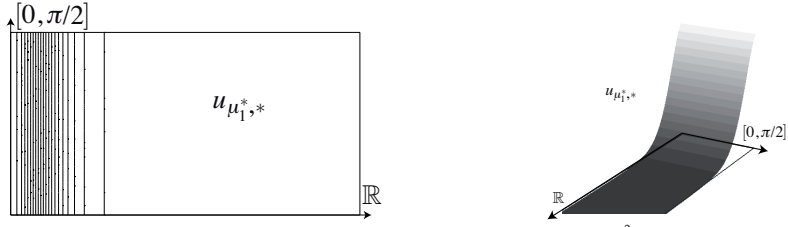}
\hspace{2cm}
\includegraphics[width=3.5cm]{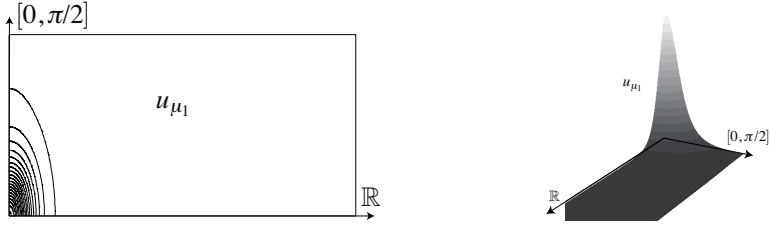}
\caption{\small\emph{Case $d=5$, $p=2.78$, $\theta=\Theta(5, 2.8)$ and $\Lambda_1=\Lambda_*^\theta(\mu_1^*)$, corresponding to the crossing of the curve $\mu\mapsto(\Lambda_*^\theta(\mu),J_*^\theta(\mu))$ with the non-symmetric curve $\mu\mapsto(\Lambda^\theta(\mu),J^\theta(\mu))$. Plot of the symmetric solution, $u_{\mu_1^*,*}$: {\rm Left.--} level lines, {\rm Right.--} 3d plot.}}
\end{center}
\end{figure}
\begin{figure}[ht]
\vspace*{-48pt}
\begin{center}
\includegraphics[width=5cm]{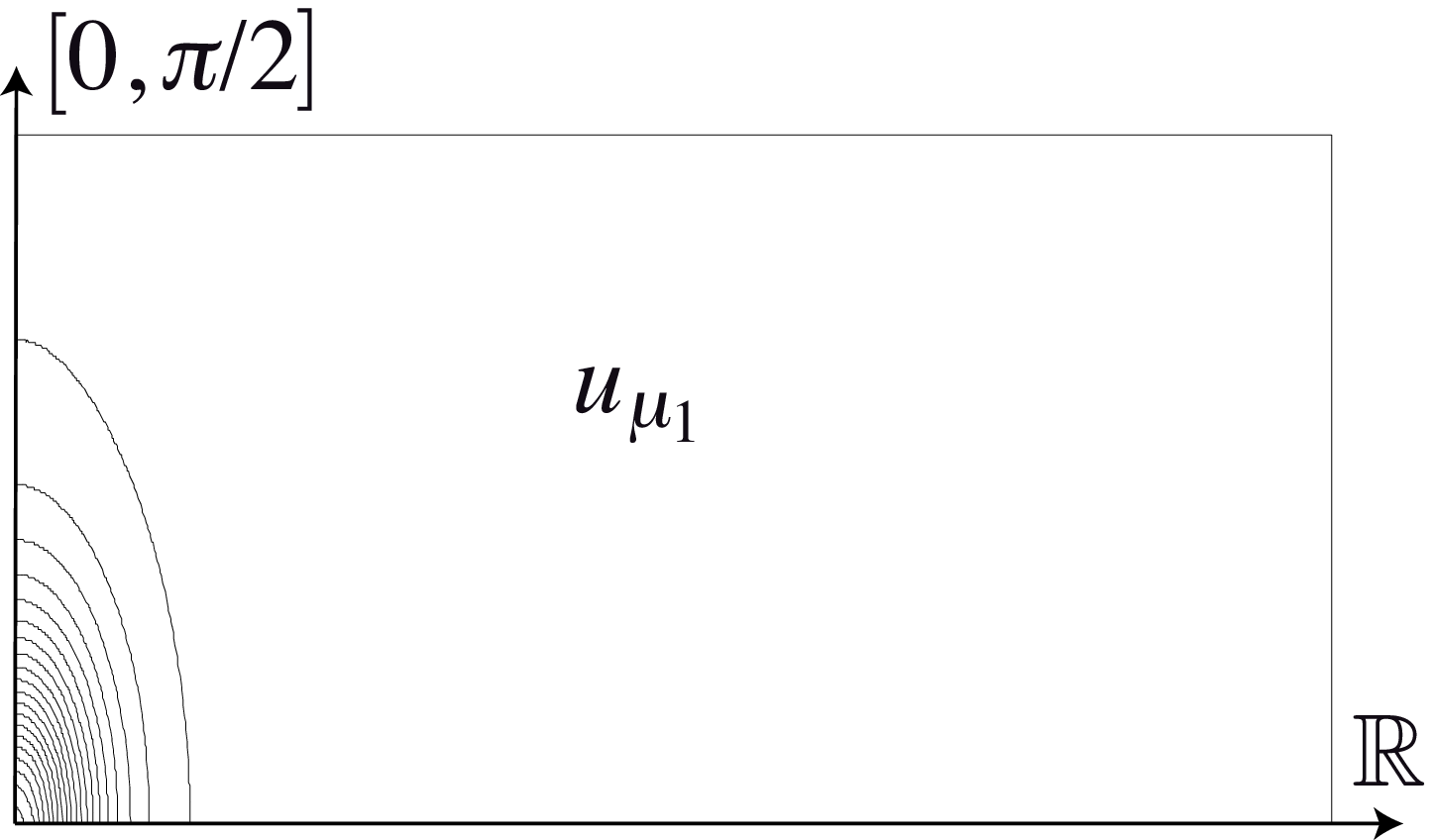}
\hspace{2cm}
\includegraphics[width=3cm]{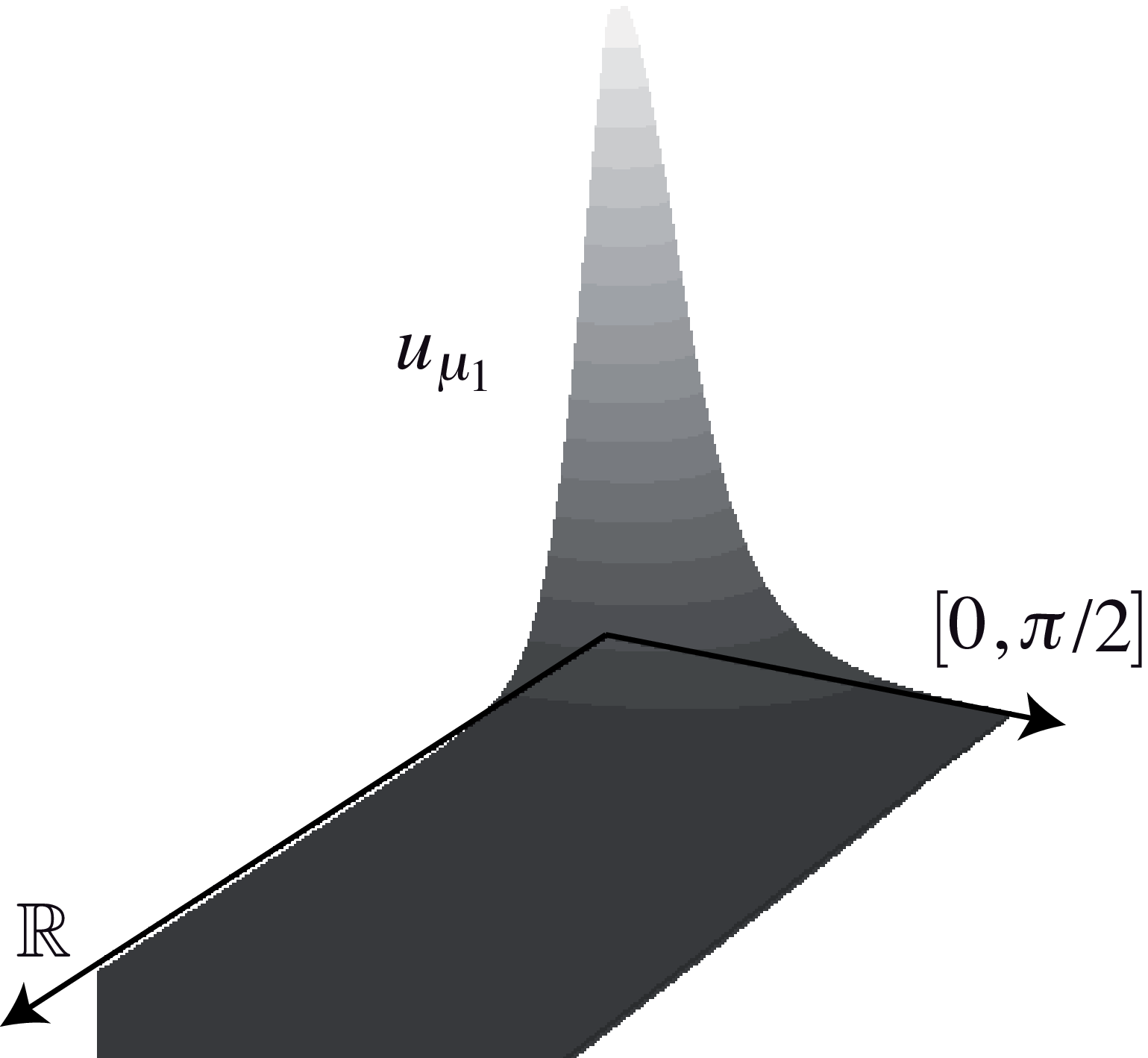}
\caption{\small\emph{Case $d=5$, $p=2.78$, $\theta=\Theta(5, 2.8)$ and $\Lambda_1=\Lambda^\theta(\mu_1)$, corresponding to the crossing of the curve $\mu\mapsto(\Lambda_*^\theta(\mu),J_*^\theta(\mu))$ with the non-symmetric curve $\mu\mapsto(\Lambda^\theta(\mu),J^\theta(\mu))$. Plot of the symmetric solution, $u_{\mu_1}$: {\rm Left.--} level lines, {\rm Right.--} 3d plot.}}
\end{center}
\vspace*{-12pt}
\end{figure}

When $p$ takes smaller values in the range $(2,p^*(\theta_0,d))$, the branch bifurcates towards the right and stays under the symmetric branch. In other words, the map $\mu\mapsto\Lambda^\theta(\mu)$ is monotone increasing for $\mu>\mu_{\rm FS}$. See Fig.~6.
\begin{figure}[hb]
\begin{center}
\includegraphics[width=5cm]{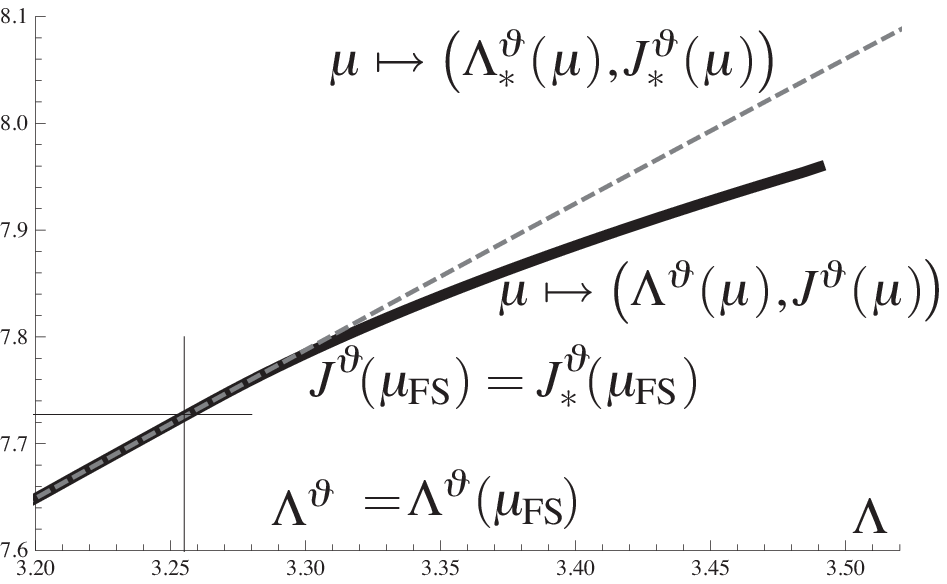}
\hspace{2cm}
\includegraphics[width=5cm]{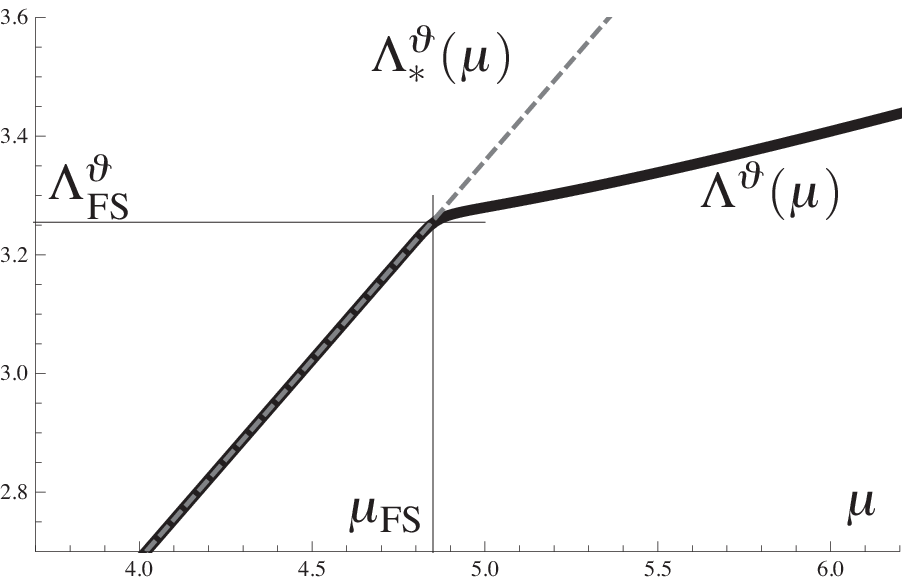}
\caption{\small\emph{Plots for $d=5$, $p=2.7$, $\theta=\Theta(5, 2.8)$. {\rm Left.--} Plot of $\mu\mapsto\(\Lambda^\theta(\mu),J^\theta(\mu)\)$ and (dashed curve) $\mu\mapsto\(\Lambda_*^\theta(\mu),J_*^\theta(\mu)\)$. {\rm Right.--} Reparametrization $\mu\mapsto\Lambda_*^\theta(\mu)$ (dashed curve) and $\mu\mapsto\Lambda^\theta(\mu)$.}}
\end{center}
\end{figure}
\clearpage
We may now come back to Fig.~1, right. An enlargement of the curves $\mu\mapsto\(\Lambda^\theta(\mu),J^\theta(\mu)\)$ clearly shows how one moves from the limiting pattern of Fig.~2 (left) to the generic case of Fig.~3 (left) and finally to the regime of Fig.~6~(left) when $\theta$ varies in $[\Theta(p,d),1)$: see Fig.~7.
\begin{figure}[ht]
\begin{center}
\includegraphics[width=5cm]{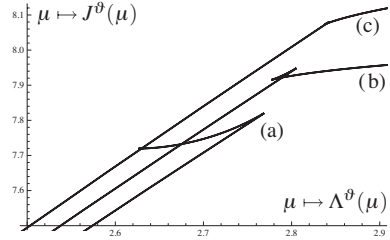}
\caption{\small\emph{Case $d=5$, $p=2.8$: curves $\mu\mapsto\(\Lambda^\theta(\mu),J^\theta(\mu)\)$ for $\theta=\Theta(2.8,5)\approx0.7143$ {\rm (a)}, $\theta\approx0.7213$ {\rm (b)} and $\theta\approx0.7283$ {\rm (c)}.}}
\end{center}
\end{figure}

\section*{Concluding remarks}

In this paper, we have observed that any critical point of $Q_\Lambda^\theta$ is also a critical point of $Q_\mu^1$ and can be rewritten as a solution of \eqref{eqmu} up to a multiplication by a constant. Using reparametrizations, it is therefore obvious that \emph{extremals} for \eqref{Ineq:CKN} belong to a union of branches that can all be parametrized by $\mu$. We have found no evidence for other branches than the ones made of symmetric solutions and of the non-symmetric ones that bifurcate from the symmetric solutions (when the number of eigenvalues of the operator $\mathcal H$ changes as $\mu$ increases). Among non-symmetric branches, the first one is the best candidate for extremal functions in \eqref{Ineq:CKN}. 

At this point we have no reason to discard the possibility of secondary bifurcations. Branches of solutions which do not bifurcate from the symmetric ones may also exist. Among the various branches, we have no theoretical reason to decide which one minimizes the energy, and the minimum may jump from one to another. This is indeed the phenomenon we have observed for instance in Fig.~3.

However, the branch that we have computed is a good candidate for minimizing the energy. It is the natural one when one starts with small values of $\mu$ and tries to optimize locally the energy functional, and it has the correct behavior as \hbox{$\mu\to\infty$}. Known estimates, like the ones of \cite{DEL2011}, show that there is not much space for unexpected solutions in the range of parameters or of the energies. It is therefore quite reasonable to conjecture that the solutions that we have computed are the actual extremals for Caffarelli-Kohn-Nirenberg inequalities and provide a complete scenario for the \emph{symmetry breaking} phenomenon of the \emph{extremals}, even if a complete proof is still missing.

\clearpage

\begin{thebibliography}{10}

\bibitem{Caffarelli-Kohn-Nirenberg-84}
Luis Caffarelli, Robert Kohn, and Louis Nirenberg.
\newblock First order interpolation inequalities with weights.
\newblock {\em Compositio Math.}, 53(3):259--275, 1984.

\bibitem{Catrina-Wang-01}
Florin Catrina and Zhi-Qiang Wang.
\newblock On the {C}affarelli-{K}ohn-{N}irenberg inequalities: sharp constants,
  existence (and nonexistence), and symmetry of extremal functions.
\newblock {\em Comm. Pure Appl. Math.}, 54(2):229--258, 2001.

\bibitem{DDFT}
Manuel del Pino, Jean Dolbeault, Stathis Filippas, and Achilles Tertikas.
\newblock A logarithmic {H}ardy inequality.
\newblock {\em J. Funct. Anal.}, 259(8):2045--2072, 2010.

\bibitem{Felli-Schneider-03}
Veronica Felli and Matthias Schneider.
\newblock Perturbation results of critical elliptic equations of
  {C}affarelli-{K}ohn-{N}irenberg type.
\newblock {\em J. Differential Equations}, 191(1):121--142, 2003.

\bibitem{MR2437030}
Jean Dolbeault, Maria~J. Esteban, and Gabriella Tarantello.
\newblock The role of {O}nofri type inequalities in the symmetry properties of
  extremals for {C}affarelli-{K}ohn-{N}irenberg inequalities, in two space
  dimensions.
\newblock {\em Ann. Sc. Norm. Super. Pisa Cl. Sci. (5)}, 7(2):313--341, 2008.

\bibitem{DEL2011}
Jean Dolbeault, Maria~J. Esteban, and Michael Loss.
\newblock Symmetry of extremals of functional inequalities via spectral
  estimates for linear operators.
\newblock {\em ArXiv e-prints, To appear in J. Math. Phys.}, September 2012.

\bibitem{DELT09}
Jean Dolbeault, Maria~J. Esteban, Michael Loss, and Gabriella Tarantello.
\newblock On the symmetry of extremals for the {C}affarelli-{K}ohn-{N}irenberg
  inequalities.
\newblock {\em Adv. Nonlinear Stud.}, 9(4):713--726, 2009.

\bibitem{1005}
Jean Dolbeault and Maria~J. Esteban.
\newblock Extremal functions for {C}affarelli-{K}ohn-{N}irenberg and
  logarithmic {H}ardy inequalities.
\newblock {\em To appear in Proc. A Edinburgh}, 2012.

\bibitem{springerlink:10.1007/s00526-011-0394-y}
Jean Dolbeault, Maria Esteban, Gabriella Tarantello, and Achilles Tertikas.
\newblock Radial symmetry and symmetry breaking for some interpolation
  inequalities.
\newblock {\em Calculus of Variations and Partial Differential Equations},
  42:461--485, 2011.

\bibitem{Oslo}
Jean Dolbeault and Maria~J. Esteban.
\newblock About existence, symmetry and symmetry breaking for extremal
  functions of some interpolation functional inequalities.
\newblock In Helge Holden and Kenneth~H. Karlsen, editors, {\em Nonlinear
  Partial Differential Equations}, volume~7 of {\em Abel Symposia}, pages
  117--130. Springer Berlin Heidelberg, 2012.
\newblock 10.1007/978-3-642-25361-4-6.

\bibitem{Landau-Lifschitz-67}
Lev~Davidovich Landau and E.~Lifschitz.
\newblock {\em Physique th\'eorique. Tome III: M\'ecanique quantique. Th\'eorie
  non relativiste. (French)}.
\newblock Deuxi\`eme \'edition. Translated from russian by E. Gloukhian.
  \'Editions Mir, Moscow, 1967.

\bibitem{DE2012}
Jean Dolbeault and Maria~J. Esteban.
\newblock Symmetry breaking in nonlinear elliptic partial differential
  equations: a scenario based on bifurcations and reparametrization.
\newblock In preparation.

\bibitem{MR0121101}
Joseph~B. Keller.
\newblock Lower bounds and isoperimetric inequalities for eigenvalues of the
  {S}chr\"odinger equation.
\newblock {\em J. Mathematical Phys.}, 2:262--266, 1961.

\bibitem{Lieb-Thirring76}
Elliott~H. Lieb and Walter Thirring.
\newblock {\em Inequalities for the moments of the eigenvalues of the
  Schr{\"o}dinger Hamiltonian and their relation to Sobolev inequalities},
  pages 269--303.
\newblock Essays in Honor of Valentine Bargmann, E. Lieb, B. Simon, A. Wightman
  Eds. Princeton University Press, 1976.

\bibitem{MR2001882}
Didier Smets and Michel Willem.
\newblock Partial symmetry and asymptotic behavior for some elliptic
  variational problems.
\newblock {\em Calc. Var. Partial Differential Equations}, 18(1):57--75, 2003.

\bibitem{Lin-Wang-04}
Chang-Shou Lin and Zhi-Qiang Wang.
\newblock Symmetry of extremal functions for the
  {C}affarelli-{K}ohn-{N}irenberg inequalities.
\newblock {\em Proc. Amer. Math. Soc.}, 132(6):1685--1691 (electronic), 2004.

\end{thebibliography}

\end{document}